\documentclass[11pt,reqno]{amsart}

\usepackage{amsmath, amsthm, amscd, amssymb, amsfonts, amsbsy}

\usepackage{esint}

\usepackage[bitstream-charter]{mathdesign}
\usepackage[T1]{fontenc}

\newtheorem{definition}{Definition}
\newtheorem{lemma}{Lemma}

\newtheorem{remark}{Remark}
\newtheorem{theorem}{Theorem}

\newcommand{\R}{\mathbb{R}}

\newcommand{\bR}{\mathbb{R}}

\newcommand{\norm}[1]{{\left\Vert #1 \right\Vert}}
\newcommand{\set}[1]{{\left\lbrace #1 \right\rbrace}}
\newcommand{\lr}[1]{{\left( #1 \right)}}

\newcommand{\intR}{\int_{\R^3}}

\newcommand{\na}{\nabla}

\newcommand{\ti}{\to \infty}

\newcommand{\ou}{\overline{u}}

\newcommand{\qand}{\quad \text{ and } \quad}
\newcommand{\qif}{\quad \text{ if } \quad}
\newcommand{\qas}{\quad \text{ as }}
\newcommand{\qfa}{\quad \text{ for all }}

\DeclareMathOperator*{\divg}{div \, }

\begin{document}
\baselineskip=18pt

\title{Logarithmic improvement of a Liouville-type theorem for the stationary Navier--Stokes equations
}

\author{Youseung Cho \& Minsuk Yang}

\address{
Yonsei University, Department of Mathematics, Yonseiro 50, Seodaemungu, Seoul 03722, Republic of Korea
}
\email{(Cho) youseung@yonsei.ac.kr \& (Yang) m.yang@yonsei.ac.kr}

\begin{abstract}
We establish a new Liouville-type theorem for the stationary Navier--Stokes equations in $\mathbb{R}^3$.
The main result is an improvement of the previous result with a logarithmic factor based on understanding of $L^p$ growth of velocity field near infinity.
\\
\noindent{\bf keywords:} 
Navier--Stokes equations, 
Liouville-type theorem, 
Energy estimates
\end{abstract}

\maketitle

\section{Introduction}
\label{S1}
%========================================

In this paper, we consider smooth solutions to the three dimensional stationary Navier--Stokes equations 
\begin{align}
\label{E11}
\begin{split}
-\Delta u + (u\cdot \na)u + \na \pi &= 0, \\
\divg u &= 0.
\end{split}
\end{align}
Here $u$ denotes the velocity of the fluid, and $\pi$ denotes the pressure.

It has been known that the solutions to \eqref{E11} are equal to zero under some additional conditions.
These results, known as Liouville-type theorems, have been actively studied within the mathematical fluid mechanics community. 
Among the many interesting research results, we only mention the ones most relevant to this paper.
Galdi \cite{MR2808162} proved a Liouville-type theorem under the assumption $u \in L^{\frac{9}{2}}(\bR^3)$.
Seregin and Wang \cite{MR3937507} proved Liouville-type theorems, assuming that the velocity field belongs to some Lorentz spaces.
In particular, $u = 0$ if there is $\frac{12}{5} < p < 3$ such that 
\[
\liminf_{R \to \infty} \frac{\norm{u}_{L^p(B(2R) \setminus B(R))}}{R^{\frac{2}{p}-\frac{1}{3}}} < \infty.
\]
Tsai \cite{MR4354995} proved Liouville-type theorems in the whole space, the half space, or a periodic slab.
In particular, $u = 0$ if there is $\frac{12}{5} \le p \le 3$ such that 
\begin{equation}
\label{E12}
\liminf_{R \to \infty} \frac{\norm{u}_{L^p(B(2R) \setminus B(R))}}{R^{\frac{2}{p}-\frac{1}{3}}} = 0.
\end{equation}
Cho et al. \cite{MR4699113} extended these results by applying an iteration method effectively.
In particular, $u = 0$ if there is $\frac{3}{2} < p < 3$ such that 
\begin{equation}
\label{E13}
\liminf_{R \to \infty} \frac{\norm{u}_{L^p(B(2R) \setminus B(R))}}{R^{\frac{2}{p}-\frac{1}{3}}} < \infty.
\end{equation}

Seregin \cite{MR3538409} proved a Liouville-type theorem based on $BMO^{-1}$ condition of the potential function.
In particular, $u = 0$ if there is a skew-symmetric tensor $V \in BMO(\R^3)$ such that $u = \divg V$.
Chae and Wolf \cite{MR3959933} proved a Liouville-type theorem based on a growth condition of the potential function over large balls, which improves the result in \cite{MR3538409}.
In particular, $u = 0$ if $u = \divg V$ and there is $3 < p \le 6$ such that 
\[
\sup_{R > 1} \frac{\norm{V-(V)_{B(R)}}_{L^p(B(R))}}{R^{\frac{2}{p}+\frac{1}{3}}} < \infty.
\]
Recently, Bang and Yang \cite{bang2024saintvenant} partially improved the assumption in Seregin \cite{MR3538409} and Chae and Wolf \cite{MR3959933}.
More precisely, $u = 0$ if there is a skew-symmetric tensor $V$ such that $u = \divg V$ and 
\[
\sup_{R > 2} \frac{\norm{V-(V)_{B(R)}}_{L^p(B(R))}}{R^{\frac{2}{p}+\frac{1}{3}}(\log R)^{{1 \over 3}-{1 \over p}}} < \infty.
\]
Motivated by this work, we improve a result in Cho et al. \cite{MR4699113} by a logarithmic factor.

Throughout the paper, we use the notation
\[
S(\rho) = \set{x \in \R^3 : \frac{3}{4} \rho \le |x| < \rho}.
\]
Here is the main result of this paper.

\begin{theorem}
\label{T1}
Let $u$ be a smooth solution to \eqref{E11}.
If there is $\frac{3}{2} < p < 3$ such that 
\begin{equation}
\label{E14}
\limsup_{\rho \to \infty} 
\frac{\norm{u}_{L^p(S(\rho))}}{\rho^{\frac{2}{p}-\frac{1}{3}} (\log \rho)^{\frac{3}{p}-1}} < \infty,
\end{equation}
then $u = 0$.
\end{theorem}

We observed that the following quantity plays a critical role in the proof of Theorem \ref{T1}:
\[
\liminf_{\rho \ti} 
\lr{\frac{\rho^{\frac{3}{p}-1}}{g(\rho)} \norm{u}_{L^3(S(\rho))}} = \infty,
\]
where $g(\rho) = \rho^{\frac{2}{p} - \frac{1}{3}} (\log \rho)^{\frac{3}{p} -1}$.
Based on this obervation, we can show that 
$1 <
\norm{u}_{L^3(S(\rho))} 
\le 
2 \norm{\ou}_{L^3(S(\rho))}$ eventually, which makes us to get an optimal result.

\begin{remark}
%\begin{enumerate}
%\item
Theorem \ref{T1} shows that if a non-trivial solution to \eqref{E11} should grow near infinity and the growth rate has improved compared to the previous results.
%\item
%If $\frac{3}{2} < p \le \frac{12}{5}$, then the proof of Theorem \ref{T1} does not rely on the iteration method. 
%\end{enumerate}
\end{remark}

Organization of paper:
We introduce general notation and preliminary lemmas in Section \ref{S2}. 
In Section \ref{S3}, we state the main lemmas of the paper related to estimating energies.
Section \ref{S4} contains the proofs of the main lemmas.
The proof of Theorem \ref{T1} using the main lemmas is then given in Section \ref{S5}.

%========================================
\section{Preliminaries}
\label{S2}
%========================================

Here is a short list of our notations.
\begin{itemize}
\item
We denote $a \lesssim b$ when there is a generic positive number $C$ such that $|a| \le C|b|$.
If an implicit constant depends on the parameter $q$, we will write $a \lesssim_q b$.
\item
If the domain of integration is the whole space $\R^3$, we often write $\int f$ instead of $\intR f(x) dx$.
\item
Let $B(r) = \set{x \in \R^3 : |x|<r}$. 
For $0 < r < R$ and $1/2 \le k \le 1$, we define $\varphi = \varphi_{r,R} \in C^1_c(B(R))$ to be a radially decreasing scalar function such that $\varphi = k$ on $B(r)$ and 
\begin{equation}
\label{E21}
(R-r)|\na \varphi| \le C < \infty,
\end{equation}
where the constant $C$ does not depend on $r$, $R$, and $k$.
\end{itemize}

\begin{definition}
Let $|E|$ denote the Lebesgue measure of a measurable set $E$ in $\bR^3$.
Define 
\[
\norm{\ou}_{L^p(E)} = \norm{u - (u)_E}_{L^p(E)}
\]
where $(u)_E = \fint_E u = \frac{1}{|E|} \int_E u$.
\end{definition}

To handle the pressure term, we shall use the following well-known lemma.

\begin{lemma}
[Lemma 1 of \cite{MR4354995}]
\label{L1}
Let $E$ be a bounded Lipschitz domain in $\bR^n$, $n \ge 2$.
Let $1 < p < \infty$ and denote $L^p_0(E) = \set{f \in L^p(E) : (f)_E = 0}$.
There is a linear map 
\[
T : L^p_0(E) \to W^{1,p}_0(E; \bR^n),
\]
such that for any $f \in L^p_0(E)$, $v = T f$ is a vector field that satisfies
\[
v \in W^{1,p}_0(E)^n, \quad \divg v = f, \quad \norm{\nabla v}_{L^p(E)} \le C(E,q) \norm{f}_{L^p(E)},
\]
where the constant $C(E,q)$ does not depend on $f$.
If $RE = \set{Rx:x \in E}$, then $C(RE,p) = C(E,p)$.
\end{lemma}

The following iteration lemma is very useful.

\begin{lemma}
[Lemma 8 of \cite{MR3720841}]
\label{L2}
Let $f(r)$ be a non-negative bounded function on $[R_0, R_1] \subset \bR_{+}$.
If there are non-negative constants $A$, $B$, $C$ and positive exponents $b < a$ and a parameter $\theta \in (0,1)$ such that for all $R_0 \le r < R \le R_1$
\[
f(r) \le \theta f(R) + \frac{A}{(R-r)^a} + \frac{B}{(R-r)^b} + C,
\]
then for all $R_0 \le r < R \le R_1$
\[
f(r) \lesssim_{a, \theta} \frac{A}{(R-r)^a} + \frac{B}{(R-r)^b} + C.
\]
\end{lemma}

%\begin{lemma}
%[Theorem 1.1 of \cite{Cho_2024}]
%\label{L3}
%If there is $\frac{3}{2} < p < 3$ such that 
%\[
%\liminf_{R \to \infty} \frac{\norm{u}_{L^p(B(2R) \setminus B(R))}}{R^{\frac{2}{p}-\frac{1}{3}}} < \infty,
%\]
%then $u = 0$.
%\end{lemma}

%========================================
\section{Energy Estimates}
\label{S3}
%========================================

In this section, we present the key lemmas related to energy estimation.
Lemma \ref{L4} provides an estimate for the derivative of the energy functional.
Lemma \ref{L5} and Lemma \ref{L5a} establish a local energy inequality. 
Lemma \ref{L6} examines the growth of energies within balls using a standard iteration lemma. 
Lemma \ref{L7} is a direct consequence of Lemma \ref{L6}.

\begin{definition}
\label{D2}
Let $\eta : [0,\infty) \to [0,1]$ be a function defined by 
\begin{align*}
\eta(t) = 
\begin{cases}
1 &\qif 0 \le t \le \frac{3}{4}, \\
-4t+4 &\qif \frac{3}{4} < t \le 1, \\
0 &\qif 1 \le t < \infty.
\end{cases} 
\end{align*}
For $\rho > 0$ define 
\[
E(\rho) = \intR |\na u(x)|^2 \eta \lr{\frac{|x|}{\rho}} dx.
\]
\end{definition}

\begin{lemma}
\label{L4}
The energy $E \in C^1((0,\infty))$ satisfies 
\[
\frac{d}{d\rho} E(\rho) 
= \frac{4}{\rho^2}
\int_{S(\rho)} 
|\na u(x)|^2 
|x| dx \qfa \rho > 0.
\]
In particular, 
\[
\frac{d}{d\rho} E(\rho) 
\ge \frac{3}{\rho} \norm{\na u}_{L^2(S(\rho))}^2 
\qfa \rho > 0.
\]
\end{lemma}

\begin{lemma}
\label{L5}
If $0 < \frac{3}{2} \rho \leq r < R \leq 2 \rho$, then
\[
\int |\na u|^2 \varphi_{r,R}
\lesssim 
\frac{1}{R-r}
\norm{\ou}_{L^2(S(R))}  
\norm{\na u}_{L^2(S(R))}
+ 
\frac{1}{R-r} 
\norm{u}_{L^3(S(R))} 
\norm{\ou}_{L^3(S(R))}^2
\]
where $\varphi_{r,R} \in C^1_c (B(R))$ is the function defined in section \ref{S2}.
\end{lemma}

\begin{lemma}
\label{L5a}
For all $\rho > 0$, 
\[
E(2\rho) \lesssim 
\rho^{-1}
\norm{\ou}_{L^2(S(2\rho))}  
\norm{\na u}_{L^2(S(2\rho))}
+ 
\rho^{-1} 
\norm{u}_{L^3(S(2\rho))} 
\norm{\ou}_{L^3(S(2\rho))}^2.
\]
\end{lemma}

We note that the implied constants in Lemma \ref{L5} and Lemma \ref{L5a} are absolute numbers.

\begin{lemma}
\label{L6}
If $\frac{3}{2} < p < 3$, then for all $\rho > 0$
\[
\norm{\na u}_{L^2(B(\sqrt3 \rho))}^2 
\lesssim 
\frac{1}{\rho^{\frac{6-p}{2p-3}}} 
\norm{u}_{L^p(S(2\rho))}^{\frac{3p}{2p-3}}
+ \frac{1}{\rho}.
\]
\end{lemma}

\begin{lemma}
\label{L7}
Let $\frac{3}{2} < p < 3$. 
If there is $M \ge 1$ and $\rho \ge 3$ such that 
\[
\norm{u}_{L^p(S(2\rho))} 
\le M \rho^{\frac{2}{p}-\frac{1}{3}} (\log \rho)^{\frac{3}{p}-1},
\]
then 
\[
\norm{\na u}_{L^2(B(\sqrt 3 \rho))}^2 
\lesssim M^{\frac{3p}{2p-3}} (\log \rho)^{\frac{3(3-p)}{2p-3}}.
\]
\end{lemma}

We note that the implied constants in Lemma \ref{L6} and Lemma \ref{L7} depend only on $p$.

%========================================
\section{Proofs of Main Lemmas}
\label{S4}
%========================================

In this section, we present the proofs of the key lemmas introduced in the previous section.

\subsection{Proof of Lemma \ref{L4}}

Fix $\rho > 0$ and let $0 < h < \frac{1}{8} \rho$.
A straightforward calculation confirms \begin{align*}
\eta\lr{\frac{|x|}{\rho+h}} 
- \eta\lr{\frac{|x|}{\rho}} 
= 
\begin{cases}
0 
&\qif 0 \le \frac{|x|}{\rho} < \frac{3}{4}, \\
4 \frac{|x|}{\rho} - 3 
&\qif \frac{3}{4} \le \frac{|x|}{\rho} < \frac{3}{4} (1 + \frac{h}{\rho}), \\
4 \frac{|x|}{\rho} - 4 \frac{|x|}{\rho + h}
&\qif \frac{3}{4} (1 + \frac{h}{\rho}) \le \frac{|x|}{\rho} < 1, \\
- 4 \frac{|x|}{\rho + h} + 4 
&\qif 1 \le \frac{|x|}{\rho} < 1 + \frac{h}{\rho},\\
0 
&\qif 1 + \frac{h}{\rho} \le \frac{|x|}{\rho} < \infty.
\end{cases} 
\end{align*}
Thus, we have
\begin{align*}
\frac{E(\rho + h) - E(\rho)}{h} 
&= \frac{1}{h} \intR |\na u(x)|^2 
\lr{\eta\left(\frac{|x|}{\rho+h}\right) 
- \eta\left(\frac{|x|}{\rho}\right)}
dx \\
&= \frac{1}{h} 
\int_{\set{\frac{3}{4} \le \frac{|x|}{\rho} < \frac{3}{4} (1 + \frac{h}{\rho})}} 
|\na u(x)|^2 
\lr{4 \frac{|x|}{\rho} - 3}
dx \\
&\quad + \frac{1}{h} 
\int_{\set{\frac{3}{4} (1 + \frac{h}{\rho}) \le \frac{|x|}{\rho} < 1}} 
|\na u(x)|^2 
\lr{4 \frac{|x|}{\rho} - 4 \frac{|x|}{\rho + h}}
dx \\
&\quad + \frac{1}{h} 
\int_{\set{1 \le \frac{|x|}{\rho} < 1 + \frac{h}{\rho}}} 
|\na u(x)|^2 
\lr{- 4 \frac{|x|}{\rho + h} + 4}
dx \\
&=: I_1(h) + I_2(h) + I_3(h).
\end{align*}
One can easily confirm $\lim_{h \to 0} I_1(h) = \lim_{h \to 0} I_3(h) = 0$ with a quick estimation.
Indeed, 
\begin{align*}
I_1(h) 
&= 
\int_{\set{\frac{3}{4} \le \frac{|x|}{\rho} < \frac{3}{4} (1 + \frac{h}{\rho})}} 
|\na u(x)|^2 
\frac{1}{h} \lr{4 \frac{|x|}{\rho} - 3}
dx \\
&\le 
\int_{\set{\frac{3}{4} \le \frac{|x|}{\rho} < \frac{3}{4} (1 + \frac{h}{\rho})}} 
|\na u(x)|^2 
\frac{3}{\rho} 
dx 
\to 0 
\qas h \to 0
\end{align*}
and 
\begin{align*}
I_3(h) 
&= 
\int_{\set{1 \le \frac{|x|}{\rho} < 1 + \frac{h}{\rho}}} 
|\na u(x)|^2 
\frac{1}{h} \lr{- 4 \frac{|x|}{\rho + h} + 4}
dx \\
&\le 
\int_{\set{1 \le \frac{|x|}{\rho} < 1 + \frac{h}{\rho}}} 
|\na u(x)|^2 
\frac{4}{\rho + h} 
dx 
\to 0 
\qas h \to 0.
\end{align*}
Finally, 
\begin{align*}
I_2(h) 
&= 
\int_{\set{\frac{3}{4} (1 + \frac{h}{\rho}) \le \frac{|x|}{\rho} < 1}} 
|\na u(x)|^2 
\frac{1}{h} \lr{4 \frac{|x|}{\rho} - 4 \frac{|x|}{\rho + h}}
dx \\
&= 
\int_{\set{\frac{3}{4} (1 + \frac{h}{\rho}) \le \frac{|x|}{\rho} < 1}} 
|\na u(x)|^2 
\lr{\frac{4|x|}{\rho (\rho + h)}}
dx \\
&\to  
\frac{4}{\rho^2} 
\int_{\set{\frac{3}{4} \le \frac{|x|}{\rho} < 1}} 
|\na u(x)|^2 
|x| dx \qas h \to 0.
\end{align*}
Since $\set{x \in \R^3 : \frac{3}{4}  \le \frac{|x|}{\rho} < 1} = S(\rho)$, we have 
\[
\lim_{h \to 0} \frac{E(\rho + h) - E(\rho)}{h} 
= \lim_{h \to 0} I_2(h) 
= 
\frac{4}{\rho^2}
\int_{S(\rho)} 
|\na u(x)|^2 
|x| dx.
\]
A straightforward calculation confirms \begin{align*}
\eta\lr{\frac{|x|}{\rho}} 
- \eta\lr{\frac{|x|}{\rho - h}} 
= 
\begin{cases}
0 
&\qif 0 \le \frac{|x|}{\rho} < \frac{3}{4}(1 - \frac{h}{\rho}), \\
4 \frac{|x|}{\rho - h} - 3 
&\qif \frac{3}{4}(1 - \frac{h}{\rho}) \le \frac{|x|}{\rho} < \frac{3}{4}, \\
4 \frac{|x|}{\rho - h} - 4 \frac{|x|}{\rho}
&\qif \frac{3}{4} \le \frac{|x|}{\rho} < 1 - \frac{h}{\rho}, \\
- 4 \frac{|x|}{\rho} + 4 
&\qif 1 - \frac{h}{\rho} \le \frac{|x|}{\rho} < 1, \\
0 
&\qif 1 \le \frac{|x|}{\rho} < \infty.
\end{cases} 
\end{align*}
Repeating the same calculation, we obtain that
\[
\lim_{h \to 0} \frac{E(\rho) - E(\rho - h)}{h} 
= \frac{4}{\rho^2}
\int_{S(\rho)} 
|\na u(x)|^2 
|x| dx.
\]
If $x \in S(\rho)$, then $|x| \ge \frac{3}{4} \rho$, and hence
\[
\frac{4}{\rho^2}
\int_{S(\rho)} 
|\na u(x)|^2 
|x| dx
\ge 
\frac{3}{\rho} \int_{S(\rho)} |\nabla u(x)|^2 dx. 
\]
\qed

\subsection{Proof of Lemma \ref{L5}}

If $0 < \frac{3}{2} \rho \le r < R \leq 2 \rho$, then $3R \le 6 \rho \le 4r$.
Since $\frac{3}{4} R \le r$, we have 
\[
\set{x \in \R^3 : r < |x| < R} 
\subset \set{x \in \R^3 : \frac{3}{4} R \le |x| < R}
= S(R).
\]
Let
\[
\ou = u - (u)_{S(R)} 
\qand w = T \left(\ou \cdot \nabla \varphi_{r, R} \right)
\quad \text{ in } S(R)
\]
where $T$ is the operator in Lemma \ref{L1} and $\varphi_{r, R}$ is the cut-off function satisfying \eqref{E21}, which is defined in Section \ref{S2}.
By Lemma \ref{L1}
\[
\nabla \cdot w = \ou \cdot \nabla \varphi_{r, R}
\]
and for each $1 < q < \infty$,
\begin{equation}
\label{E41}
\norm{\nabla w}_{L^q(S(R))} 
\lesssim_q 
\frac{1}{R-r} \norm{\ou}_{L^q(S(R))}.
\end{equation}
Multiplying $\ou \varphi_{r,R} - w$ to the first equation of \eqref{E11} and integrating on $\R^3$, we get 
\[
- \int \Delta u \cdot (\ou \varphi_{r,R} - w) + \int (u \cdot \na) u \cdot (\ou \varphi_{r,R} - w) = 0,
\]
where the pressure term is absent due to $\divg (\ou \varphi_{r,R} - w) = 0$.
Integration by parts and $\divg u = 0$ yield
\[
\int |\na u|^2 \varphi_{r, R} 
= 
-\int \na u :(\ou \otimes \na \varphi_{r,R})
+ \int \na u : \nabla w 
+ \frac{1}{2} \int | \ou|^2 u \cdot \nabla \varphi_{r, R}
- \int (u \otimes \ou) : \nabla w,
\]
where $F:G = \sum_{i,j} F_{ij} G_{ij}$ for tensor-valued function $F$ and $G$.
Notice that the supports of $\nabla \varphi_{r, R}$ and $w$ are contained in $S(R)$.
Using the H\"older inequality and \eqref{E41} with $q=2$, we get 
\begin{align*}
&-\int \na u :(\ou \otimes \na \varphi_{r,R})
+ \int \na u : \nabla w \\
&\lesssim 
\frac{1}{R-r} \norm{\ou}_{L^2(S(R))} \norm{\na u}_{L^2(S(R))}
+ \norm{\na u}_{L^2(S(R))} 
\norm{\nabla w}_{L^2(S(R))} \\
&\lesssim 
\frac{1}{R-r} \norm{\ou}_{L^2(S(R))}
\norm{\na u}_{L^2(S(R))}.
\end{align*}
Using the H\"older inequality and \eqref{E41} with $q = 3$, we get
\begin{align*}
&\frac{1}{2} \int | \ou|^2 u \cdot \nabla \varphi_{r,\overline R}
- \int (u \otimes \ou) : \nabla w \\
&\lesssim 
\frac{1}{R-r} \norm{\ou}_{L^{3}(S(R))}^2 
\norm{u}_{L^{3}(S(R))}
+ \norm{u}_{L^{3}(S(R))} 
\norm{\ou}_{L^{3}(S(R))} 
\norm{\nabla w}_{L^{3}(S(R))} 
\\
&\lesssim 
\frac{1}{R-r} \norm{u}_{L^{3}(S(R))}
\norm{\ou}_{L^{3}(S(R))}^2.
\end{align*}
Combining these two estimates with the above identity gives the result.
\qed

\subsection{Proof of Lemma \ref{L5a}}

Let $n > 8$ and 
\begin{align*}
\eta_n(t)
= 
\begin{cases}
1 - \frac{4}{n}
&\qif 0 \le t < \frac{3}{4}, \\
-2n \left( t - \frac{3}{4} \right)^2 + 1 - \frac{4}{n} 
&\qif \frac{3}{4} \le t < \frac{3}{4} + \frac{1}{n}, \\
-4t + 4 - \frac{2}{n}
&\qif \frac{3}{4} + \frac{1}{n} \le t < 1 - \frac{1}{n}, \\
2n(t-1)^2 
&\qif 1 - \frac{1}{n} \le t < 1,\\
0 
&\qif 1 \le t < \infty.
\end{cases} 
\end{align*}
Then 
\begin{align*}
\eta'_n(t)
= 
\begin{cases}
0
&\qif 0 \le t < \frac{3}{4}, \\
-4n \left( t - \frac{3}{4} \right)
&\qif \frac{3}{4} \le t < \frac{3}{4} + \frac{1}{n}, \\
-4
&\qif \frac{3}{4} + \frac{1}{n} \le t < 1 - \frac{1}{n}, \\
4n(t-1)
&\qif 1 - \frac{1}{n} \le t < 1,\\
0 
&\qif 1 \le t < \infty.
\end{cases} 
\end{align*}
Note that
\[
\eta_n(t) \nearrow \eta(t) \le 1 \qand -\eta_n'(t) \nearrow -\eta'(t) \le 4
\]
where $\eta(t)$ is the function defined in Definition \ref{D2}.
Since $\eta_n \in C^1_c ((0,1))$, we may use Lemma \ref{L5} with $r = 3\rho/2$, $R = 2 \rho$, and $\varphi_{3\rho/2, 2\rho} = \eta_n(|x|/2\rho)$ to get
\[
\int |\na u|^2 \eta_n \left( \frac{|x|}{2\rho} \right)
\lesssim 
\rho^{-1}
\norm{\ou}_{L^2(S(2\rho))}  
\norm{\na u}_{L^2(S(2\rho))}
+ 
\rho^{-1}
\norm{u}_{L^3(S(2\rho))} 
\norm{\ou}_{L^3(S(2\rho))}^2.
\]
We get the desired result by observing
\[
\int |\na u|^2 \eta_n \left( \frac{|x|}{2\rho} \right)
\nearrow 
\int |\na u|^2 \eta \left( \frac{|x|}{2\rho} \right) 
= E(2\rho).
\]
\qed

\subsection{Proof of Lemma \ref{L6}}

Let $\frac{3}{2} < p < 3$ and $0 < \sqrt 3 \rho \leq r < R \leq 2 \rho$.
Denote 
\[
A(R) = \set{x \in \R^3 : \frac{\sqrt 3}{2} R \le |x| < R}.
\]
We notice that $A(R) \subset S(2\rho)$.
Modifying the support of the test function and corresponding averages in the proof of Lemma \ref{L5}, one can easily obtain
\begin{align*}
\norm{\na u}_{L^2(B(r))}^2 
&\le 
C \frac{1}{R-r} 
\norm{\ou}_{L^2(A(R))} 
\norm{\na u}_{L^2(A(R))}
+
C \frac{1}{R-r} 
\norm{u}_{L^3(A(R))} 
\norm{\ou}_{L^3(A(R))}^2 \\
&=: J_1 + J_2.
\end{align*}
By the H\"older inequality 
\begin{align*}
J_1 
&\lesssim 
\frac{1}{R-r} 
|A(R)|^{1/6}
\norm{\ou}_{L^3(A(R))} 
\norm{\na u}_{L^2(A(R))} \\
&\lesssim 
\lr{\frac{1}{(R-r)^{4}} \rho^3}^{\frac{1}{6}}
\lr{\frac{1}{R-r} \norm{\ou}_{L^3(A(R))}^3}^{\frac{1}{3}}
\norm{\na u}_{L^2(A(R))}.
\end{align*}
By the Young inequality 
\[
J_1 
\le 
\frac{1}{6} 
\norm{\na u}_{L^2(A(R))}^2
+ 
C \frac{1}{R-r} 
\norm{\ou}_{L^3(A(R))}^3
+
C \frac{1}{(R-r)^{4}} \rho^3.
\]
Since 
\[
\norm{u}_{L^3(A(R))} \le \norm{\ou}_{L^3(A(R))} 
+ |A(R)|^{{1 \over 3}-{1 \over p}} \norm{u}_{L^p(A(R))},
\]
we have 
\[
J_2 
\lesssim 
\frac{1}{R-r} 
\norm{\ou}_{L^3(A(R))}^3
+
\frac{1}{R-r} 
\frac{1}{\rho^{{3-p \over p}}}
\norm{u}_{L^p(A(R))}
\norm{\ou}_{L^3(A(R))}^2.
\]
Combining the estimates for $J_1$ and $J_2$, we get 
\begin{equation}
\label{E42}
\begin{split}
\norm{\na u}_{L^2(B(r))}^2 
&\le 
\frac{1}{6} \norm{\na u}_{L^2(A(R))}^2
+ 
C \frac{1}{R-r} 
\norm{\ou}_{L^3(A(R))}^{3}  \\ 
&\quad +  
C \frac{1}{R-r} 
\frac{1}{\rho^{{3-p \over p}}}
\norm{u}_{L^p(A(R))}
\norm{\ou}_{L^3(A(R))}^2 
+
C \frac{1}{(R-r)^{4}} 
\rho^3.
\end{split}
\end{equation}
By the H\"older inequality, the Poincar\'e--Sobolev inequality, $A(R) \subset S(2\rho)$, and $A(R) \subset B(R)$, we get
\[
\norm{\ou}_{L^3(A(R))} 
\le 
\norm{\ou}_{L^p(A(R))}^{p \over 6-p} 
\norm{\ou}_{L^6(A(R))}^{6-2p \over 6-p} 
\lesssim 
\norm{\ou}_{L^p(A(R))}^{p \over 6-p} 
\norm{\na u}_{L^2(A(R))}^{6-2p \over 6-p} 
\lesssim 
\norm{u}_{L^p(S(2\rho))}^{p \over 6-p} 
\norm{\na u}_{L^2(B(R))}^{6-2p \over 6-p}.
\]
Thus, \eqref{E42} becomes 
\begin{align*}
\norm{\na u}_{L^2(B(r))}^2 
&\le 
\frac{1}{6} \norm{\na u}_{L^2(B(R))}^2
+ 
C \frac{1}{R-r}
\norm{u}_{L^p(S(2\rho))}^{3p \over 6-p} 
\norm{\na u}_{L^2(B(R))}^{3(6-2p) \over 6-p} \\
&\quad + 
C \frac{1}{R-r} 
\frac{1}{\rho^{{3-p \over p}}}
\norm{u}_{L^p(S(2\rho))}^{6+p \over 6-p} 
\norm{\na u}_{L^2(B(R))}^{2(6-2p) \over 6-p}
+
C \frac{1}{(R-r)^{4}}
\rho^3.
\end{align*}
Since $\frac{3}{2} < p < 3$, we can use the Young inequality to get 
\begin{equation}
\label{E43}
\begin{split}
\norm{\na u}_{L^2(B(r))}^2 
&\le 
\frac{1}{2} \norm{\na u}_{L^2(B(R))}^2
+ 
C \frac{1}{(R-r)^{6-p \over 2p-3}} 
\norm{u}_{L^p(S(2\rho))}^{3p \over 2p-3} \\
&\quad + 
C \frac{1}{(R-r)^{6-p \over p}} 
\frac{1}{\rho^{{(3-p)(6-p) \over p^2}}}
\norm{u}_{L^p(S(2\rho))}^{6+p \over p}
+
C \frac{1}{(R-r)^{3}}
\rho^4.
\end{split}
\end{equation}
Since $0 < \frac{(2p-3)(6+p)}{3p^2} < 1$ and $R-r < \rho$, we use the Young inequality to get 
\begin{align*}
&\frac{1}{(R-r)^{6-p \over p}} 
\frac{1}{\rho^{{(3-p)(6-p) \over p^2}}}
\norm{u}_{L^p(S(2\rho))}^{6+p \over p} \\
&= 
\lr{\frac{1}{(R-r)^{6-p \over 2p-3}} 
\norm{u}_{L^p(S(2\rho))}^{3p \over 2p-3}}^{\frac{(2p-3)(6+p)}{3p^2}}
\lr{\frac{(R-r)^2}{\rho^3}}^{1-\frac{(2p-3)(6+p)}{3p^2}} \\
&\lesssim 
\frac{1}{(R-r)^{6-p \over 2p-3}} 
\norm{u}_{L^p(S(2\rho))}^{3p \over 2p-3}
+ \frac{1}{\rho}.
\end{align*}
Thus, \eqref{E43} becomes 
\begin{equation}
\label{E44}
\norm{\na u}_{L^2(B(r))}^2 
\le 
\frac{1}{2} \norm{\na u}_{L^2(B(R))}^2
+ 
C \frac{1}{(R-r)^{6-p \over 2p-3}} 
\norm{u}_{L^p(S(2\rho))}^{3p \over 2p-3} 
+
C \frac{1}{(R-r)^{4}}
\rho^3.
\end{equation}
By Lemma \ref{L2}, we obtain that for $0 < \sqrt 3 \rho \leq r < R \leq 2 \rho$
\[
\norm{\na u}_{L^2(B(r))}^2 
\lesssim 
\frac{1}{(R-r)^{6-p \over 2p-3}} 
\norm{u}_{L^p(S(2\rho))}^{3p \over 2p-3} 
+
\frac{1}{(R-r)^{4}}
\rho^3.
\]
Taking $r = \sqrt 3 \rho$ and $R = 2 \rho$, we get the desired result.
\qed

\subsection{Proof of Lemma \ref{L7}}

Since
\[
\norm{u}_{L^p(S(2\rho))} 
\le M \rho^{\frac{2}{p}-\frac{1}{3}} (\log \rho)^{\frac{3}{p}-1},
\]
we obtain from Lemma \ref{L6} that 
\begin{align*}
\norm{\na u}_{L^2(B(\sqrt3 \rho))}^2 
&\lesssim 
\frac{1}{\rho^{\frac{6-p}{2p-3}}} 
\norm{u}_{L^p(S(2\rho))}^{\frac{3p}{2p-3}}
+ \frac{1}{\rho} \\
&\le 
\frac{1}{\rho^{\frac{6-p}{2p-3}}} 
(M \rho^{6-p \over 3p} (\log \rho)^{3-p \over p})^{3p \over 2p-3} 
+ 1 \\
&\le 
2 M^{\frac{3p}{2p-3}} 
(\log \rho)^{\frac{9-3p}{2p-3}}.
\end{align*}
\qed

%========================================
\section{Proof of Theorem \ref{T1}}
\label{S5}
%========================================

\noindent
\textbf{Step 1.}

\noindent
Let $3/2 < p < 3$ and
\begin{equation}
\label{E50}
g(\rho) = \rho^{\frac{2}{p} - \frac{1}{3}} (\log \rho)^{\frac{3}{p} -1}.
\end{equation}
Then the condition \eqref{E14} becomes 
\begin{equation}
\label{E51}
\limsup_{\rho \ti} \frac{\norm{u}_{L^p(S(\rho))}}{g(\rho)} < \infty.
\end{equation}
Note that 
\[
\rho^{-\frac{1}{3}}
\norm{u}_{L^3(S(\rho))} 
=
\lr{\rho^{\frac{2}{3}-\frac{3}{p}} g(\rho)}
\lr{\frac{\rho^{\frac{3}{p}-1}}{g(\rho)} \norm{u}_{L^3(S(\rho))}}
\]
and 
\[
\lim_{\rho \ti} \rho^{\frac{2}{3}-\frac{3}{p}} g(\rho) 
= \lim_{\rho \ti} \frac{(\log \rho)^{\frac{3}{p} -1}}{\rho^{\frac{1}{p}-\frac{1}{3}}}
= 0.
\]
Thus, $\liminf_{\rho \ti} \lr{\frac{\rho^{\frac{3}{p}-1}}{g(\rho)} \norm{u}_{L^3(S(\rho))}} < \infty$ implies 
\[
\liminf_{\rho \ti}  
\rho^{-\frac{1}{3}} 
\norm{u}_{L^3(S(\rho))} 
= 0,
\]
which gives $u = 0$ thanks to \eqref{E12} (Theorem 1 of Tsai \cite{MR4354995}).
This is the desired goal.
Therefore, we only need to consider the case 
\begin{equation}
\label{E52}
\liminf_{\rho \ti} 
\lr{\frac{\rho^{\frac{3}{p}-1}}{g(\rho)} \norm{u}_{L^3(S(\rho))}} = \infty.
\end{equation}

\noindent
\textbf{Step 2.}

\noindent
From now, we assume \eqref{E52}.
The condition \eqref{E51} implies that there are $M > 1$ and $R_1 > 3$ such that 
\begin{equation}
\label{E53}
\frac{\norm{u}_{L^p(S(\rho))}}{g(\rho)} \le M \qfa \rho \ge R_1.
\end{equation}
The condition \eqref{E52} implies that there exist $R_2 > R_1$ such that 
\begin{equation}
\label{E54}
2M < \frac{\rho^{\frac{3}{p}-1}}{g(\rho)} \norm{u}_{L^3(S(\rho))} \qfa \rho \ge R_2.
\end{equation}
Since 
\[
\lim_{\rho \ti} \rho^{\frac{2}{3} - \frac{1}{p}} (\log \rho)^{\frac{3}{p} -1} = \infty,
\]
there is $R_3 > R_2$ such that 
\[
1< 
2M \rho^{\frac{2}{3} - \frac{1}{p}} (\log \rho)^{\frac{3}{p} -1}
= 2 M \frac{g(\rho)}{\rho^{\frac{3}{p}-1}}
\qfa \rho \ge R_3.
\]
Thus, from \eqref{E54}, we have 
\[
1< 
2 M \frac{g(\rho)}{\rho^{\frac{3}{p}-1}}  
< \norm{u}_{L^3(S(\rho))} 
\qfa \rho \ge R_3.
\]
Using \eqref{E53} and \eqref{E54}, we have for all $\rho \ge R_3$
\begin{align*}
\norm{u}_{L^3(S(\rho))} 
&\le 
\norm{\ou}_{L^3(S(\rho))} 
+ 
|S(\rho)|^{{1 \over 3} - {1 \over p}}
\norm{u}_{L^p(S(\rho))} \\
&\le 
\norm{\ou}_{L^3(S(\rho))} 
+ 
\frac{M g(\rho)}{\rho^{\frac{3}{p}-1}} \\
&< 
\norm{\ou}_{L^3(S(\rho))} 
+ 
{1 \over 2} \norm{u}_{L^3(S(\rho))}.
\end{align*}
Therefore, 
\begin{equation}
\label{E55}
1 <
\norm{u}_{L^3(S(\rho))} 
\le 
2 \norm{\ou}_{L^3(S(\rho))} \qfa \rho \ge R_3.
\end{equation}

\noindent
\textbf{Step 3.}

\noindent
By Lemma \ref{L5a} 
\begin{equation}
\label{E56}
\begin{split}
E(2\rho) 
&\lesssim 
\rho^{-1} 
\norm{\ou}_{L^2(S(2\rho))} 
\norm{\na u}_{L^2(S(2\rho))} 
+
\rho^{-1} 
\norm{\ou}_{L^3(S(2\rho))}^2
\norm{u}_{L^3(S(2\rho))} \\
&=: K_1 + K_2.
\end{split}
\end{equation}
We claim that
\begin{equation}
\label{E57}
K_1 
\lesssim \left( \log \rho \right)^{\frac{9-3p}{6-p}} \norm{\na u}_{L^2(S(2\rho))}^{\frac{18-6p}{6-p}}.
\end{equation}
If $3/2 < p \le 2$,  we use \eqref{E55}, the Poincar\'e--Sobolev inequality, and \eqref{E53} to get
\begin{align*}
K_1 
&\lesssim \rho^{-1} \norm{\ou}_{L^2(S(2\rho))} 
\norm{\na u}_{L^2(S(2\rho))} \norm{\ou}_{L^3(S(2\rho))} \\
&\lesssim \rho^{-1} \norm{u}_{L^p(S(2\rho))}^{3p \over 6-p} 
\norm{\na u}_{L^2(S(2\rho))}^{18-6p \over 6-p} \\
&\lesssim \left( \log \rho \right)^{\frac{9-3p}{6-p}}
\norm{\na u}_{L^2(S(2\rho))}^{18-6p \over 6-p}.
\end{align*}
If $2 < p \le 12/5$, we use the Poincar\'e inequality, the Jensen inequality, and \eqref{E53} to get
\begin{align*}
K_1 
&\lesssim \rho^{-\frac{4p-6}{6-p}} \norm{\ou}_{L^2(S(2\rho))}^{\frac{4p-6}{6-p}} 
\norm{\na u}_{L^2(S(2\rho))}^{\frac{18-6p}{6-p}} \\
&\lesssim \left( \rho^{-\frac{3}{p} + \frac{1}{2}} \norm{u}_{L^p(S(2\rho))} \right)^{\frac{4p-6}{6-p}} 
\norm{\na u}_{L^2(S(2\rho))}^{\frac{18-6p}{6-p}} \\
&\lesssim \rho^{-\frac{2p-3}{3p}} \left( \log \rho \right)^{\frac{(4p-6)(3-p)}{p(6-p)}} \norm{\na u}_{L^2(S(2\rho))}^{\frac{18-6p}{6-p}} \\
&\lesssim \norm{\na u}_{L^2(S(2\rho))}^{\frac{18-6p}{6-p}}.
\end{align*}
If $12/5 < p < 3$, we have $\frac{18-6p}{6-p} < 1$. 
So we use the Jensen inequality, Lemma \ref{L7}, and \eqref{E53} to get
\begin{align*}
K_1 
&\lesssim \rho^{-\frac{p}{3} + \frac{1}{2}} 
\norm{u}_{L^p(S(2\rho))} 
\norm{\na u}_{L^2(S(2\rho))}^{\frac{18-6p}{6-p}}
\left( \log \rho \right)^{\frac{3(5p-12)(3-p)}{(6-p)(2p-3)}} \\
&\lesssim \rho^{-\frac{1}{p} + \frac{1}{6}} 
\left( \log \rho \right)^{\frac{3-p}{p} + \frac{3(5p-12)(3-p)}{(6-p)(2p-3)}} 
\norm{\na u}_{L^2(S(2\rho))}^{\frac{18-6p}{6-p}} \\
&\lesssim \norm{\na u}_{L^2(S(2\rho))}^{\frac{18-6p}{6-p}}.
\end{align*}
In any cases, \eqref{E57} is true.
From Lemma \ref{L4} 
\begin{align*}
K_1 
&\lesssim \left( \log \rho \right)^{\frac{9-3p}{6-p}} \norm{\na u}_{L^2(S(2\rho))}^{\frac{18-6p}{6-p}} \\
&\lesssim \left( \rho \log \rho E'(2\rho) \right)^{9-3p \over 6-p}. 
\end{align*}
We use \eqref{E55} and the Poincar\'e--Sobolev inequality, \eqref{E53} to get
\[
K_2 
\lesssim 
\rho^{-1} 
\norm{\ou}_{L^3(S(2\rho))}^3 
\lesssim 
\rho^{-1} 
\norm{u}_{L^p(S(2\rho))}^{3p \over 6-p} 
\norm{\na u}_{L^2(S(2\rho))}^{18-6p \over 6-p}.
\]
From Lemma \ref{L4} 
\begin{align*}
K_2 
&\lesssim \rho^{-1} 
g(2\rho)^{3p \over 6-p} 
\left( \rho E'(2\rho) \right)^{9-3p \over 6-p} \\
&\lesssim \left( \rho \log \rho E'(2\rho) \right)^{9-3p \over 6-p}.
\end{align*}
Therefore, 
\[
E(\rho) 
\lesssim \left( \rho \log \rho E'(\rho) \right)^{9-3p \over 6-p} \qfa \rho \ge R_3.
\]

\noindent
\textbf{Step 4.}

\noindent
Note that $E(\rho)$ is increasing.
Suppose there is $R_4 > R_3$ such that $E(\rho) > 0$ for all $\rho \ge R_4$.
Then,
\[
\frac{1}{\rho \log \rho}
\lesssim 
\frac{E'(\rho)}{E(\rho)^{{6-p \over 9-3p}}}
\qfa \rho \ge R_4.
\]
Integrating both sides, we get 
\[
\int_{R_4}^{R_5}
\frac{1}{\rho \log \rho}
d\rho
\lesssim 
\int_{R_4}^{R_5}
\frac{E'(\rho)}{E(\rho)^{{6-p \over 9-3p}}}
d\rho
\lesssim 
\frac{1}{E(R_4)^{{2p-3 \over 9-3p}}}
-
\frac{1}{E(R_5)^{{2p-3 \over 9-3p}}}.
\]
Letting $R_5 \to \infty$ gives a contradiction. 
Since $E(\rho)$ is increasing, $E(\rho) = 0$ for all $\rho$, which implies $\na u = 0$.
Due to \eqref{E14}, we get $u=0$.
This completes the proof of Theorem \ref{T1}.
\qed

%========================================
%\section*{Acknowledgement}
%
%This work was supported by the National Research Foundation of Korea(NRF) grant funded by the Korean government(MSIT) (No. 2021R1A2C4002840).

\section*{Declarations}

\textbf{Conflict of interest}
The authors declare that there is no conflict of interest regarding the publication of this article.

\textbf{Ethics approval}
The work does not involve the study of live subjects.
Ethical approval is not needed for the study.

\textbf{Funding}
The work is supported by the National Research Foundation of Korea(NRF) grant funded by the Korean government(MSIT) (No. 2021R1A2C4002840).

\textbf{Data availability} 
Data sharing not applicable to this article as no datasets were generated or
analysed during the current study.

\end{document}